\let\accentvec\vec
\documentclass{llncs}

\let\vec\accentvec
\usepackage{makeidx}  

\usepackage{amsmath}
\usepackage{dsfont}
\usepackage{verbatim}
\usepackage{graphicx}
\usepackage{color}
\usepackage{amssymb}

\usepackage{url}
\usepackage{hyperref}

\newcommand{\Ecal}{\mathcal{E}}

\newcommand{\Gcal}{\mathcal{G}}
\newcommand{\Mcal}{\mathcal{M}}
\newcommand{\Ncal}{\mathcal{N}}
\newcommand{\Pcal}{\mathcal{P}}
\newcommand{\Xcal}{\mathcal{X}}
\newcommand{\Ycal}{\mathcal{Y}}
\newcommand{\Mixt}{\operatorname{Mixt}}
\newcommand{\supp}{\operatorname{supp}}

\newcommand{\RBM}{\ensuremath{\operatorname{RBM}}}

\newcommand{\DBN}{\ensuremath{\operatorname{DBN}}}
\newcommand{\MPD}{\ensuremath{\Mcal}}    
\newcommand{\UMPD}{\ensuremath{\Mcal^{\ast}}} 

\newcommand{\ol}{\overline}

\begin{document}
\title{Maximal Information Divergence  from Statistical Models defined by Neural Networks} 
\titlerunning{Maximal Divergence from Neural Networks}  
%
\author{Guido Mont\'ufar\inst{1} \and Johannes Rauh\inst{2} \and Nihat Ay\inst{2,3}}
%
%
\tocauthor{Guido Mont\'ufar, Johannes Rauh, and Nihat Ay}
\institute{Department of Mathematics, Pennsylvania State University, \\University Park, PA 16802, USA,\\
\email{gfm10@psu.edu},\\ 
\and
Max Planck Institute for Mathematics in the Sciences,\\
Inselstra\ss e 22, 04103 Leipzig, Germany,\\
\email{$\{$jrauh,nay$\}$@mis.mpg.de},\\
\and
Santa Fe Institute, 1399 Hyde Park Road, Santa Fe, NM 87501, USA.
}

\maketitle              

\begin{abstract}
We review recent results about the maximal values of the Kullback-Leibler information divergence from statistical models defined by neural networks, including na\"ive Bayes models, restricted Boltzmann machines, deep belief networks, and various classes of exponential families. 
We illustrate approaches to compute the maximal divergence from a given model starting from simple sub- or super-models. 
We give a new result for deep and narrow belief networks with finite-valued units. 
\keywords{neural network, exponential family, Kullback-Leibler divergence, multi-information}
\end{abstract}

\section{Introduction}

In statistical learning theory, probability models are used to infer representations of data. 
In model selection it is often assumed that the model approximation errors are negligible compared with the statistical approximation errors.
This assumption may not always be justified in practice; in some cases even full dimensional models only fill a small portion of the space of probability distributions, and telling the general structure of the data generating distributions, in order to constrain the possible model classes, is difficult. 

\looseness-1 
Here we take a complementary perspective, disregarding the statistical approximation errors and focussing on the model
approximation errors. We quantify the model approximation error of a model $\Mcal$ by the divergence function $p\mapsto D(p\|\Mcal)=\inf_{q\in\Mcal}D(p\|q)$, 
where $D(p\|q)=\sum_xp(x)\log\frac{p(x)}{q(x)}$  is the Kullback-Leibler divergence from $p$ to $q$.\footnote{
We formulate our results in such a way that they are independent from the logarithm's base used in the definition of the  divergence. 
 } 
We study the maximum value of $D(\cdot\|\Mcal)$, which corresponds to a worst-case analysis. 
The ideas from this paper can also be used to study the expectation value given a prior on the set of target distributions, see~\cite{wupes2012}. 
The model approximation error can be used as a criterion for model selection. 
Related ideas are discussed in~\cite{AyMontufarRauhICCN2011} in the context of model design and reinforcement learning.

Most probability models with hidden variables are singular and not identifiable. 
Moreover, data distributions that are not contained in these models can have several maximum likelihood estimates. 
Although controlling parameter-identifiability 
is crucial when estimating learning coefficients in Bayesian model selection, we will instead focus on the value of the data likelihood and the sets of maximizing distributions, irrespective of their parameters. 

In general, the function $D(\cdot\|\Mcal)$ has no explicit formula, making the estimation of the
maximizers and the maximum value difficult. 
For exponential families the situation is slightly better, 
as for each distribution $p$ the divergence $D(p\|\cdot)$ has a unique minimizer
over~$\overline{\Mcal}$.  For certain families, such as independence models and convex exponential families,
there even is a closed formula for this function.  The approximation properties of various classes of exponential
families have been studied in~\cite{Matus2003,matusiid,AyKnauf06:Maximizing_Multiinformation,Rauh11:Finding_Maximizers,Rauh13:Optimal_Expfams,wupes2012,Juricek10:Maximization_from_multinomial_distributions}. 
The divergence from complicated models can be estimated by finding tractable exponential subfamilies.
This idea was used in~\cite{NIPS2011_0307} to study approximation errors of restricted Boltzmann machines. 

The representational power of neural networks has been studied for many years and by too many authors to refer to appropriately at this place, see for instance~\cite{Cybenko1988b,DBLP:journals/nn/HornikSW89,Funahashi1998209}.  
The representational power of the networks discussed in this paper has been studied, in particular,  in~\cite{LeRoux2008,Hinton2008,LeRoux2010,Montufar2011,NIPS2011_0307,Montufar2010a}. 

Section~\ref{sec:maxim-inform-diverg} reviews bounds on $D_{\Mcal}$ for statistical models defined by
neural networks and for exponential families.  Section~\ref{sec:estim-inform-diverg} discusses strategies to
bound~$D_{\Mcal}$ via sub-models and super-models, and discusses a class of exponential families contained in restricted Boltzmann machines and deep belief networks. Section~\ref{sec:discussion} puts our results in perspective.

\section{Maximal information divergence} 
\label{sec:maxim-inform-diverg}

\begin{figure}[t]
\setlength{\unitlength}{\textwidth}
\centering
\begin{picture}(.8,.26)(0,-.035)
\put(0,-.05){\includegraphics[trim=3.7cm 20.3cm 4.7cm 2cm, clip=true, width=\unitlength]{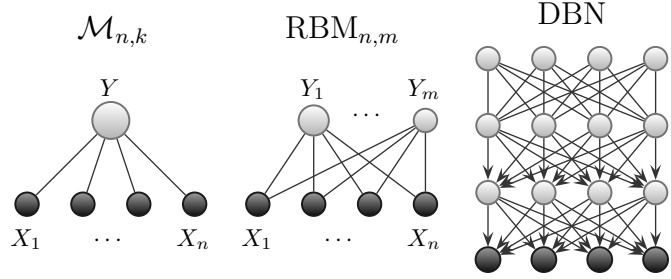}}
\put(.025,.21){\begin{minipage}{.2\unitlength}\center {\large$\Mcal_{n,k}$}\end{minipage}} 

\put(.275,.21){\begin{minipage}{.2\unitlength}\center {\large$\RBM_{n,m}$}\end{minipage}}

\put(.525,.23){\begin{minipage}{.2\unitlength}\center {\large$\DBN$}\end{minipage}}
\end{picture}
\caption{The na\"ive Bayes model~$\Mcal_{n,k}$, the restricted Boltzmann machine $\RBM_{n,m}$, and a deep belief network. 
Light (dark) nodes represent hidden (visible) variables.}
\label{graphs}
\end{figure}

We consider neural networks with a set of visible units~$X_{1},\dots,X_{n}$, where each $X_{i}$ takes values in a finite set~$\Xcal_{i}$ of cardinality $|\Xcal_i|=N_i$. See Fig.~\ref{graphs}. 
The visible state space of such a system is $\Xcal=\Xcal_{1}\times\dots\times\Xcal_{n}$.
For any subset $A\subseteq[n]$ let $N_{A}=\prod_{i\in A}N_i$ be the number of joint states of the units indexed by~$A$, and let~$N=N_{[n]}=|\Xcal|$. 
We denote the set of all probability distributions $(p_x)_{x\in\Xcal}$ on $\Xcal$ by~$\Delta(\Xcal)$, or~$\Delta$ if $\Xcal$ is understood. 
The {\bf\em maximal information divergence} from a model $\Mcal\subseteq\Delta(\Xcal)$ is $D_\Mcal:=\max_{p\in\Delta}D(p\|\Mcal)$. 
An {\em $rI$-projection} of $p\in\Delta$ onto $\Mcal$ 
is a point $p_{\Mcal}$ in the closure $\overline{\Mcal}$ of $\Mcal$ with
$D(p\|\Mcal) = D(p\|p_{\Mcal})$. 

\subsection{Probability models defined by neural networks}
\looseness-1 
The {\bf\em independence model} $\Ecal^1_n$ of $n$ variables $X_1,\ldots,X_n$ is the set of probability distributions of the form $p(x)=\prod_{i\in[n]} p_i(x_i)$ for all $x=(x_{1},\dots,x_{n})\in\Xcal$. 
This model describes non-interacting stochastic variables. 
The following result is due to Ay and Knauf~\cite[Corollary~4.10]{AyKnauf06:Maximizing_Multiinformation}. 

\begin{lemma}
\label{lemmadeayknauf}
The maximal divergence to $\Ecal^1_n$ is bounded by 
\begin{equation*}
D_{\Ecal^1_n}\leq \log(N/\max_{i\in[n]}{N_i})\;. 
\end{equation*}
If all variables are $q$-ary, then $D_{\Ecal^1_n}=(n-1)\log(q)$, and the maximizers are the uniform distributions on $q$-ary codes of cardinality $q$ and minimum distance $n$. 
\end{lemma}

The {\bf\em mixture of product distributions}~$\Mcal_{n,k}$, or {\bf\em na\"ive Bayes model}, is the graphical model on
a star graph, where the leaves are visible variables, and the internal node is a hidden variable with $k$ states. 
\begin{theorem}
  \label{thm:max-KL-from-Mix}
  Let $A\subseteq[n]$. If $k\geq N_{[n]\setminus A}$, then $D_{\Mcal_{n,k}}$ is bounded by 
  \begin{equation*}
  D_{\Mcal_{n,k}}\leq \log ( N_{A} / \max_{j\in A}N_j )\;.
  \end{equation*}
When all visible variables are binary, we have the tighter bound 
  \begin{equation*}
    D_{\Mcal_{n,k}} \leq 
  \big(   n -\left\lfloor \log_2(k) \right\rfloor - \frac{k}{2^{\left\lfloor\log_2 (k)\right\rfloor}} \big) \log(2) \;. 
    \label{boundlog}
  \end{equation*}
\end{theorem}
Note the similarity of the bounds given in Lemma~\ref{lemmadeayknauf} and Theorem~\ref{thm:max-KL-from-Mix},
In fact Theorem~\ref{thm:max-KL-from-Mix} can be derived from Lemma~\ref{lemmadeayknauf}, together with Lemma~\ref{lem:rI-of-mixture} given below. 

\medskip 

The {\bf\em restricted Boltzmann machine} $\RBM_{n,m}$ is the undirected stochastic network with full bipartite interaction graph $K_{n,m}$, where an independent set of $m$ units is hidden, and an independent set of $n$ units is visible. 
\begin{theorem} 
 \label{thm:max-KL-RBM-qary}
  \label{thm:max-KL-RBM}
Let $A\subseteq[n]$, and let $M_{1},\dots,M_{m}$ be the sizes of the state spaces of the hidden variables. 
If $1+\sum_{j\in[m]}(M_{j}-1) \geq N_{[n]\setminus A}$, then 
\begin{equation*}
D_{\RBM_{n,m}} \leq \log ( N_{A} / \max_{j\in A}N_j )\;.
\end{equation*}
When all units are binary, and $m\leq 2^{n-1}-1$, we have the tighter bound 
  \begin{equation*}
    D_{\RBM_{n,m}}\leq 
\big(     n -\left\lfloor \log_2(m+1) \right\rfloor - \frac{m+1}{2^{\left\lfloor\log_2 (m+1)\right\rfloor}}  \big) \log(2) \;.\label{boundlog1}
  \end{equation*}
  \end{theorem}

Theorem~\ref{thm:max-KL-RBM-qary} subsumes divergence bounds for na\"ive Bayes models (when $m=1$) and independence models (when $m=0$). 
This result was shown in the binary case in~\cite[Theorem~2]{NIPS2011_0307} and in the non-binary case in~\cite[Theorem~29]{montufar2013discrete}. 

\medskip
A {\bf\em deep belief network} (DBN) is a layered stochastic network with undirected bipartite interactions between the units in the deepest two layers, which form an RBM, and directed bipartite interactions between all other pairs of subsequent layers, directed towards the first layer, which is the only visible layer.

\begin{theorem} 
\label{thm:max-KL-DBN}
Consider a DBN with $L$ layers, each layer containing $n$ units with state spaces of cardinalities $q_1,\ldots,q_n$. Let $m$ be any integer with $\prod_{j=m+2}^{n}q_j\leq m\leq n$, and let $q_1\geq \cdots\geq q_m$.  
If $L\geq 2+\frac{q_1^S-1}{q_1-1}$ for some $S\in\{0,1,\ldots,m\}$, then 
\begin{equation*}
D_{\DBN}\leq \log(N_{[m-S]})\;. 
\end{equation*}

\medskip
\noindent
 In particular, when all units are binary and the network has $L\geq 1+2^S$ layers of size $n =2^{k-1}+k$, for some $S\in\{0,1, \ldots, 2^{k-1}\}$, 
 then 
 \begin{equation*}
 D_{\DBN} \leq \big( 2^{k-1} - S  \big)\log(2) \;. 
 \end{equation*}
\end{theorem}

The binary case is~\cite[Theorem~2]{Montufar2011}, together with~\cite[Theorem~18]{Kernels2012}. 
The non-binary case is new (details in~\cite{FiniteDBNs}). 

\medskip 

The bounds in Theorems~\ref{thm:max-KL-from-Mix}, \ref{thm:max-KL-RBM} and~\ref{thm:max-KL-DBN} vanish when the number of hidden units is large enough (depending on their state spaces). In this case, the models can approximate all probability distributions on the states of their visible units arbitrarily well, i.e., they are {\em universal approximators}.

All these theorems can be proved using the same strategy: First, a family of exponential sub-models is identified,
and then, the divergence from the union of these sub-models is bounded from above, as in
Theorem~\ref{thm:max-KL-from-prod-mix} below. 

\subsection{Exponential families}
Exponential families are widely-used statistical models. Examples include log-linear models, hierarchical models, and independence models. 
The information divergence maximization problem is by far better understood for exponential families than for other probability models. 
We use exponential families to approximate probability models with hidden variables. 

\medskip 

Let $\varrho=\{A_{1},\dots,A_{m}\}$ be a partition of $\Xcal$.  The {\bf\em partition model} $\Pcal_{\varrho}$ \label{partitionmodelpg} consists
of all $p\in\Delta$ with $p(x)=p(y)$ whenever $x,y$ belong to the same block of~$\varrho$; that is, the conditional distribution of~$p$, conditioned on any block~$A_{i}\in\varrho$, equals the uniform distribution. 
Partition models are the convex exponential families that contain the uniform distribution.  The following is a special case of~\cite[Corollary~1]{Matus2003}:
\begin{lemma} 
  \label{lem:part-mod-max-KL}  
  Let $\varrho=\{A_{1},\dots,A_{m}\}$ be a partition of $\Xcal$, and denote by $c(\varrho) = \max_{i\in[m]}|A_{i}|$ the coarseness of~$\varrho$.  Then $D_{\Pcal_{\varrho}} = \log(c(\varrho))$,  and the global maximizers are the distributions $p$ with $\supp(p)\cap A_i\leq 1$ for all $i\in[m]$, where $\supp(p)\cap A_i= 1$ holds only if $|A_i|=c(\varrho)$.
\end{lemma}

The maximal divergence from any exponential family of dimension $k$ can be bounded from below as follows, see~\cite[Theorem~28]{Rauh13:Optimal_Expfams}: 
\begin{theorem} 
\label{optifam}
  Let $\Ecal$ be an exponential family of dimension~$k$.  Then 
  \begin{equation*}
  D_{\Ecal}\geq\log(N) -\log(k+1)\;. 
  \end{equation*}
   If equality holds, then
  $\Ecal$ is a partition model with homogeneous partition. 
\end{theorem}

Probability models defined as marginals of exponential families can behave very different from proper exponential families.
Any finite subset of $\Delta$ can be embedded in a projection of a two-dimensional  exponential family, see~\cite{AyMontufarRauhICCN2011}:
\begin{lemma} 
\label{lemmaiccn}
Given any finite set  of probability distributions $\{p^{(i)}\}_{i=1}^K\subset\Delta_{N-1}$, 
there is a two-dimensional exponential family $\Ecal\subseteq\Delta_{K-1}$, and a linear map $\psi\colon \Delta_{K-1}\to\Delta_{N-1}$, such that $\psi(\overline{\Ecal})\supseteq\{p^{(i)}\}_{i=1}^K$. 
\end{lemma}

\section{Estimating the information divergence}
\label{sec:estim-inform-diverg}
\subsection{Subfamilies and superfamilies}

If $\Mcal'\subseteq\Mcal$ then $D_{\Mcal}\le D_{\Mcal'}$.  In special cases it is possible to have equality.

\begin{lemma}
  \label{lemma:sub}
  If $\Mcal'\subseteq\Mcal$ and if $p$ is a maximizer of the divergence from~$\Mcal$ such that $\Mcal'$ contains an
  $rI$-projection $p_{\Mcal}$ of $p$ to~$\Mcal$, then $p$ maximizes the divergence from $\Mcal'$ among the set
  $\{q\in\Delta: q_{\Mcal}\in\Mcal'\text{ for some $rI$-projection }q_{\Mcal}\}$. 
\end{lemma}
The lemma is useful for exponential families; where the set of distributions whose $rI$-project to $\Mcal$ lies in $\Mcal'$, can be parametrized via $\Mcal'+\Ncal$, where $\Ncal$ is the normal space of~$\Mcal$.
The following argument due to Jur\'i\v{c}ek~\cite{Juricek10:Maximization_from_multinomial_distributions} is an example:


Let $\Mcal=\Ecal^{1}_{n}$ be the independence model of $n$ $q$-ary variables and let $\Mcal'$ be the set of
i.i.d.~distributions.  By Lemma~\ref{lemmadeayknauf}, the uniform distribution $p$ on the states
$(1,\dots,1),(2,\dots,2),\dots,(q,\dots,q)$ maximizes the divergence from~$\Mcal$, and it is exchangeable.  Since the $rI$-projections of the set of exchangeable
distributions to $\Mcal$ belong to $\Mcal'$, Lemma~\ref{lemma:sub} implies that $p$ maximizes the divergence from
$\Mcal'$ among the exchangeable distributions, with divergence $D(p\|\Mcal') = (n-1)\log(q)$.  Now, $\Mcal'$ as a subset
of the exchangeable simplex can be identified with the multinomial model.  This proves the following
result~\cite[Theorem~1.1]{Juricek10:Maximization_from_multinomial_distributions}.
\begin{theorem}
  The maximal divergence from the multinomial model of $n$ $q$-ary variables is equal to~$(n-1)\log(q)$.
\end{theorem}



Conversely, simple subfamilies can be used to study larger models:

\begin{lemma}\label{lem:maxdivusingsub-models}
Let $\Ecal$ be an exponential family. Let $\Mcal_i$ be a sub-model of $\Ecal$ with $D_{\Mcal_i}=K$ and divergence maximizers $\Gcal_i$, for all $i\in[k]$. 
If there is a point $p\in\Gcal=\cap_i\Gcal_i$ with $p_\Ecal\in\cup_i\Mcal_i$, 
then $D_\Ecal=K$ and the divergence maximizers are exactly the points in $\Gcal$ whose $rI$-projections onto $\Ecal$ lie in $\cap_i\Mcal_i$. 
\end{lemma}

Lemma~\ref{lem:maxdivusingsub-models} can be used to prove the homogeneous case of Lemma~\ref{lemmadeayknauf} as follows: 
The independence model of $n$ $q$-ary variables contains the partition model $\Pcal_{i}$ with partition blocks $\{x\colon x_i=y_i\}$ for all $y_i\in\Xcal_i$, for any $i\in[n]$. 
By Lemma~\ref{lem:part-mod-max-KL}, the maximal divergence from the partition model $\Pcal_i$ is $D_{\Pcal_i}=(n-1)\log(q)$, and the set of maximizers is the set $\Gcal_i$ of distributions $p$ whose support $\supp(p)=\{x^{(j)}\}_j$ satisfies $x^{(j)}_i\neq x^{(j')}_i$ for all $j\neq j'$. 
The intersection $\Gcal=\cap_i\Gcal_i$ is the set of probability distributions with support on a code of minimum distance $n$. 
The $rI$-projection of an arbitrary element $p\in\Gcal$ lies in $\cap_i\Pcal_i=\{u\}$ if and only if $p$ is a uniform
distribution on a code of minimum distance $n$ and cardinality $q$.
By Lemma~\ref{lem:maxdivusingsub-models} these are the global divergence maximizers from $\Ecal_n^1$.

\subsection{Mixtures of exponential families with disjoint supports}

The \emph{mixture} $\Mixt(\Mcal_{1},\dots,\Mcal_{k})$ of $k$ models $\Mcal_{1},\dots,\Mcal_{k}\subseteq\Delta$ is the set
of probability distributions of the form $p = \sum_{i=1}^{k}\lambda_{i}p^{(i)}$, where $\lambda\in\Delta_{k-1}$ and
$p^{(i)}\in\Mcal_{i}$ for all $i\in[k]$. 
In general, mixtures are difficult to describe, even for simple models $\Mcal_{1},\dots,\Mcal_{k}$. 
The situation is much simpler when mixing models supported on disjoint subsets of~$\Xcal$: 
\begin{lemma}
  \label{lem:rI-of-mixture}
  Let $\{A_{1},\dots,A_{k}\}$ be a partition of $\Xcal$ and let $\Mcal_{1},\dots,\Mcal_{k}$ be statistical models with $\Mcal_{i}\subseteq\Delta(A_{i})$.  
  For any $p\in\Delta(\Xcal)$, the $rI$-projections of $p$ to $\Mixt(\Mcal_{1},\dots,\Mcal_{k})$ are the distributions of the form 
  \begin{equation*}
    p_{\Mcal}(x) = p(A_{i})p_{\Mcal_i}(x),\qquad\text{for all $x\in A_{i}$}\quad\text{for all $i\in[k]$},
  \end{equation*}
where $p_{\Mcal_i}$ denotes an $rI$-projection of $p(x|A_{i})$ to $\Mcal_{i}$ for all $i\in[k]$. 
\end{lemma}

We call a set $\Ycal\subseteq\Xcal_1\times\cdots\times\Xcal_n$ {\em cubical} if it can be written as a product $\Ycal=\Ycal_1\times\cdots\times\Ycal_n$ with $\Ycal_i\subseteq\Xcal_i$ for all $i\in[n]$. 
A set $\Ycal$ is cubical iff there exists a product distribution $p$ with $\supp(p) := \{x\in\Xcal:p(x)>0\} = \Ycal$ (in this
case $\Ycal_{i}=\supp(p_{i})$).  We call a partition \emph{cubical} if it consists of cubical blocks.  For any cubical
set $\Ycal$ let $\ol{\Ecal^{1}_{\Ycal}}$ denote the set of product distributions with support~$\Ycal$.

\medskip

Let $\varrho=\{A_1,\ldots,A_{k}\}$ be a cubical partition of $\Xcal$. The {\bf\em mixture of products with disjoint supports} $\boldsymbol{\varrho}$ is the model $\MPD_{\varrho}:=\Mixt(\ol{\Ecal^1_{A_1}},\ldots,\ol{\Ecal^1_{A_k}})\subseteq\Mcal_{n,k}$. 
For this kind of models, Lemmas~\ref{lemmadeayknauf} and~\ref{lem:rI-of-mixture} show: 
\begin{corollary}
  \label{cor:ind-mix-max-KL}
  Let $\varrho=\{A_1,\ldots, A_k\}$ be a cubical partition of $\Xcal$ with blocks $A_i=\Ycal_{i,1}\times\cdots\times\Ycal_{i,n}$ with $|\Ycal_{i,j}|\in\{1,q_i\}$ for all $j\in[n]$, for all $i\in[k]$. Then
  \begin{equation*}
      D_{\MPD_{\varrho}}=  \max_{i\in[k]} \log(|A_{i}|/q_i) \;.
    \end{equation*}
\end{corollary}

\subsection{Unions of exponential families}

Let $\UMPD_{n,k}=\bigcup_{\varrho:|\varrho|=k}\MPD_{\varrho}\subseteq\Mcal_{n,k}$ be the union of mixtures of products with disjoint supports $\varrho$, where $\varrho$ runs over all cubical partitions of $\Xcal$ with $k$ blocks. The set $\UMPD_{n,k}$ is not an exponential family, but a finite union of exponential families.
Similarly, let $\Mcal_{n,k,0}^{\ast}=\bigcup_{\varrho:|\varrho|=k}\Pcal_\varrho$ be the union of all partition models $\Pcal_{\varrho}$ of partitions $\varrho$
with $k$ cubical blocks. 

Our motivation for studying unions of mixture models and unions of partition models comes from the following two
results. For simplicity, we consider binary units; analogue results for non-binary units can be found in~\cite{montufar2013discrete} and~\cite{FiniteDBNs}. 

\begin{theorem}[{\cite[Theorem~1]{NIPS2011_0307}}]
\label{thm:UMPD-in-RBM}
 The binary model $\RBM_{n,m}$ contains any mixture of one arbitrary product distribution, $m-k$ product distributions with mutually disjoint supports, 
and $k$ distributions with support on any edges of the $n$-cube, for any $0\leq k\leq m$.  In particular, 
$\RBM_{n,m}$ contains $\UMPD_{n,m+1}$. 
\end{theorem}

\begin{theorem}[{\cite[Theorem~17]{Kernels2012}}]
\label{thm:PM-in-DBN}
Let $L\in\mathbb{N}$, let $k$ be the largest integer for which $L\geq 1+ 2^{(2^{k-1})}$, and let $K= 2^{k-1} + k \leq n$. 
The binary deep belief network model  with $L$ layers of width $n$ contains 
any partition model $\Pcal_\varrho$ with partition $\varrho=\{\{x\colon x_\lambda=y_\lambda \} \colon {y_\lambda\in\{0,1\}^K}\}$, where $\lambda\subseteq[n], |\lambda|=K$. 
\end{theorem}

Unions of exponential families are more difficult to describe than exponential families, but the maximal  $rI$-projection can be approximated as follows: 

\begin{theorem} 
\label{thm:max-KL-from-prod-mix}
  Let $\Xcal=\{0,1\}^n$. If $k\le 2^{n-1}$, then 
\vspace{-.2cm}  
  \begin{equation*}
    D_{\UMPD_{n,k}}\le \big(  n -\left\lfloor \log_2(k) \right\rfloor - \frac{k}{2^{\left\lfloor\log_2 (k)\right\rfloor}} \big) \log(2) \;.
  \end{equation*}
\vspace{-.2cm}  
If $k\le 2^n$, then 
  \begin{equation*}
    D_{\Mcal_{n,k,0}^\ast}\le 
    \big(  n+1 -\left\lfloor \log_2(k) \right\rfloor - \frac{k}{2^{\left\lfloor\log_2 (k)\right\rfloor}} \big) \log(2) \;.
  \end{equation*}
\end{theorem}
The first part was shown in~{\cite[Theorem~2]{NIPS2011_0307}}. 
The second part can be proved with a direct adaptation of the same proof. 
Theorem~\ref{thm:max-KL-from-prod-mix}, together with Theorems~\ref{thm:UMPD-in-RBM} and~\ref{thm:PM-in-DBN}, proves the `tighter bounds' in 
Theorems~\ref{thm:max-KL-from-Mix} and~\ref{thm:max-KL-RBM}.

\section{Discussion}
\label{sec:discussion}
When we plot the approximation error bounds of the model classes discussed here against the corresponding number of model parameters, we find that 
they all behave similarly; they all decay logarithmically on a large scale. 
This is the optimal maximal approximation error behaviour of exponential families (Theorem~\ref{optifam}). 
The bounds for partition models, homogeneous independence models, and mixtures of products with disjoint homogeneous supports, are tight. 
The na\"ive Bayes model bound is tight for many choices of the $N_i$ in the sense that it vanishes iff the model is a universal approximator, see~\cite{Montufar2010a}. 
The other bounds for the more complicated models are probably not tight. 
It is reasonable to expect that fixing the number of parameters, models with many hidden units fill the probability simplex more evenly than their counterparts with fewer or no hidden units (see, e.g.,  Lemma~\ref{lemmaiccn}). 
For the discussed model classes, this paper does not give conclusive answers in that direction, since the only maximal divergence lower-bounds are for exponential families. 
It should be mentioned, however, that the mere existence of universal approximators within a given class of networks is not always obvious and sometimes false. For example, DBNs with too narrow hidden layers are never universal approximators, regardless of their parameter count.

\appendix

\medskip
\noindent
{\bf Acknowledgement:} J.~R.\ is supported in part by the VW Foundation; G.~M.\ by DARPA grant FA8650-11-1-7145.  
\bibliography{referenzen}{}\bibliographystyle{abbrv}

\end{document}